\documentstyle [12pt]{article}

\newtheorem{theorem}{Theorem}[section]
\newtheorem{prop}[theorem]{Theorem}

\newtheorem{coroll}[theorem]{Corollary}

\newtheorem{lemm}[theorem]{Lemma}

\newtheorem{anoprop}[theorem]{Proposition}

\newenvironment{mproof}{\begin{trivlist}\item[]{\em
Proof: }}{\hfill$\Box$\end{trivlist}}

\newtheorem{eg}{\rm\sl \uppercase{Example}}[section]

\newcommand{\co}{{\cal O}}

\def \IR{\hbox{{\rm I}\kern-.2em\hbox{{\rm R}}}}
\def \iR{\hbox{{\sevenrm I\kern-.2em\hbox{\sevenrm R}}}}
\def \IN{\hbox{{\rm I}\kern-.2em\hbox{\rm N}}}
\def \IC{\hbox{{\rm I}\kern-.6em\hbox{\bf C}}}
\def \IQ{\hbox{{\rm I}\kern-.6em\hbox{\bf Q}}}
\def \ZZ{\hbox{{\rm Z}\kern-.4em\hbox{\rm Z}}}

\newcommand{\ga}{{\alpha}}

\newcommand{\ot}{\otimes}

\newcommand{\pr}{\prime}

\begin{document}

\begin{center} {\Large \bf Lexicographic Semigroupoids}

\vspace{.3in}

{\Large S. C. Power}
\vspace{.15in}

\rm
\it Department of Mathematics and Statistics\\
Lancaster University\\
England LA1 4YF
\rm
\end{center}
\vspace{.3in}
\begin{center} \bf Abstract
\end{center}
{\small
The natural lexicographic semigroupoids associated with
Cantor product spaces indexed by  countable linear orders
are classified. Applications are given to the classification
of triangular operator algebras which are direct limits of upper
triangular matrix algebras. }
\rm

\vspace{.3in}
Consider a Cantor space which is presented explicitly as an infinite
product of finite topological spaces. The
product presentation provides an equivalence relation $R$ consisting
of the pairs $(x,y)$ of points $x$ and $y$ which disagree in at most
finitely many coordinates.
This equivalence relation supports a natural locally compact totally
disconnected topology which makes  $R$ a principal groupoid.
It is well-known that in the case of countable products
such topological equivalence relations are classified by the generalised
integer obtained from the formal product of the cardinalities of the
component spaces. Furthermore, this classification is closely related
to the classification of C*-algebras that are
infinite tensor products of matrix algebras, the so-called UHF
C*-algebras. See, for example, Renault \cite{ren} and Power \cite{scp-book}.

In the present paper we consider antisymmetric topological binary
relations which are the lexicographic topological subrelations
arising from
infinite products indexed
by  general countable linear orderings. These natural  semigroupoids
are classified and their automorphism groups determined. This and related
results enable us
to give applications to the classification of triangular
operator algebras which are themselves
lexicographic products in an algebraic sense.

The  binary relations
may also be viewed as the (semigroupoid) lexicographic
products of total orderings on
finite sets, and in fact our methods are applicable
to lexicographic products of
connected antisymmetric finite partial orders.
Although applications to approximately finite
operator algebras gives our primary motivation
it seems clear that lexicographic subrelations are
interesting in their own right.

In the preliminary section  we recall how the generalised integer
associated with the presentation of the Cantor space gives
a complete invariant for the associated approximately finite groupoid.
In section 2 we classify the lexicographic products in the case of
indexing by a countable dense order. It is interesting that the
proof here is much more elementary than the case of indexing by $\ZZ$
which is taken up in section 3. Indeed the classification in the former case
is less subtle and suggests, a posteriori, that the component coordinates
must be accessible in purely order-topological terms. The proof in section 3
follows the order-topological methods
in Power \cite{scp-outer} where
automorphisms of the associated triangular algebras were studied.
These algebras
- the so-called alternation algebras - have been
considered by a number of authors,
namely Poon \cite{poo}, Hopenwasser and Power \cite{hop-scp}, and most
recently,
as part of a wider analysis, by Donsig and Hopenwasser \cite{don-hop}.
Also in section 3 we identify the automorphism group of a lexicographic
product over $\ZZ$.
In section 4 we obtain the classification for the case of
general countable linear orders and in the last section we give applications
to operator algebras.
\newpage

\section{ Preliminaries}

Let $\mu$  be a positive integer-valued function defined on a countable
linear ordering $\Omega$.
For notational convenience we denote
the discrete finite topological
space $\{1, \dots, n\}$ by [$n$].
Associate with $\mu$  the Cantor space

\[ X_{\mu} = \prod_{w \in \Omega} [r_w]
\]

\noindent  where $r_w = \mu (w)$. Write
$\tilde{R}_{\mu}$ for the equivalence relation described in the
introduction. Let $R_{\mu}$ be the antisymmetric subrelation of
points
$(x, y) \in \tilde{R}_{\mu}$,
for which $x$ preceeds $y$ in the lexicographic ordering. Thus $(y,x) \in
R_{\mu}$ if and only if $x  = (x_w)$ and $y  = (y_w)$,
the coordinates $(x_w)$  and $ (y_w)$ are equal except in at most a finite
number of coordinates, and $y_w < x_w$ for the first
index $w$ where $x$ and $y$ differ.

The basic open-closed sets for the Cantor space $X_{\mu}$ are provided by the
cylinder sets arising from the specification of a finite number of
coordinates. The topology on $\tilde{R}_{\mu}$ arises in the following
similar way. Let $F \subseteq \Omega$ be a finite subset, let
$x^{\prime}_w, y^{\prime}_w$ be specified coordinates for $w \in F$,
and let

\[
E = \{(x,y) : x_w = y_w \hbox{ for } w \notin F, x_w =
x^{\prime}_w \hbox{ and } y_w = y^{\prime}_w \hbox{ for } w \in F\}.
\]

\noindent The totality of these sets gives a base for the groupoid topology
on $\tilde{R}_\mu$ and the subrelation $R_\mu$ carries the relative
topology.

Notice that if $\pi_l, \pi_r : \tilde{R}_\mu \rightarrow X_\mu$ are the
natural coordinate projections then $\pi_l$ and $\pi_r$ are continuous
and are one to one when restricted to a basic open-closed
set $E$. General open-closed sets with this property are called {\it $G$-sets}
 and
these  are important in the following two ways.
Firstly they provide a class of sets which are invariant under {\it groupoid
isomorphism},
that is, a binary relation isomorphism that is also a homeomorphism.
Secondly they are used in the formulation of invariant
measures on the underlying space.
As a consequence groupoid isomorphisms conserve the invariant measures.

An {\it invariant measure} for a
principal groupoid $\tilde{R}$ is a Borel measure $\lambda$ for the
underlying topological space such that $\lambda
(\pi_l (G)) = \lambda (\pi_r (G))$ for every $G$-set $G$. It is a simple
matter to check that $\tilde{R}_\mu$ possesses a unique invariant
probability measure, namely the infinite product of the probability
measures $\lambda_\omega$ on $[r_\omega]$ which assign equal mass
$r_\omega^{-1}$ to each singleton set. At this stage we can deduce that if
$\tilde{R}_\mu$ is isomorphic to
the topological equivalence relation
$\tilde{R_\nu}$, associated with
$\nu : \Lambda \to \IN$, then the generalised
integers for $\mu$ and $\nu$ agree. Indeed, the hypothesised
isomorphism is a bijection $\alpha: X_\mu \rightarrow X_\nu$ such that
the map $\alpha^{(2)} : \tilde{R}_\mu \to \tilde{R}_\nu$ is a bijection and a
topological isomorphism. Since the invariant probability measures for
$\tilde{R}_\mu$ and $\tilde{R}_\nu$ are unique they must correspond
under $\alpha$, and from this it follows that they have the same range
on open-closed sets.
Thus the set of rationals $k/r$, with $k \in \ZZ$ and
$r = r_{w_1}r_{w_2}\dots r_{w_n}$ for
some $n \in \IN$, coincides with the corresponding set of rationals for $\nu$.
Equivalently

\[
 \prod_{w \in \Omega}\mu (w) \ = \
 \prod_{w \in \Lambda}\nu (w)
\]

\noindent as generalised integers.

In the antisymmetric context we shall use the argument above in a local
way (with various restrictions of $R_\mu$ in place of $\tilde{R}_\mu$)
to obtain local comparisons of the data for two given lexicographic
binary relations.

For $x$ in $X_\mu$ the (one-sided) {\it orbit} of $x$ is the set $\co (x)$
of points $y$ with $(y,x) \in R_\mu$, and the {\it closed orbit} of $x$ is
$\overline{\co (x)}$ the closure of this set. Note that if $x =
(x_w)$ with $x_w = 1$ for all $w \leq w_1$ and $x_w = \mu (w)$ for all
$w > w_1$, for some $w_1$ in $\Omega$, then $\overline{\co (x)}$ is
the set of points $y = (y_w)$ with $y_w =1$ for $w \leq w_1$.
This will be indicated by setting $y_{w_1} = 1$ and
writing

\[
\overline{\co (x)} = (..1..{|y_{w_1}|}..*..).
\]

A similar such shorthand is used in the next section to indicate
basic open-closed sets in $X_\mu$.
Note that for the particular point $x$ above the restriction $R_\mu
\mid \overline{\co (x)}$ is  isomorphic to a lexicographic ordering.

\section{ Countable dense orders}

In this section we classify the topological binary relations $R_\mu$ in
the cases when $\Omega$ is a dense linear ordering and $\mu : \Omega
\rightarrow \IN$ satisfies $\mu (w) \geq 2$ for all $w$ in $\Omega$.
There are only four such linear orderings and these
correspond to the presence and absence of maximal and minimal elements.
\vspace{.3in}

\noindent {\bf Theorem 2} \it Let $\mu : \Omega \rightarrow \IN
\backslash \{1\}$ and
$\nu : \Omega^\prime \rightarrow \IN \backslash \{1\}$ be functions on
the countable dense linear orderings  $\Omega$ and
$\Omega^\prime$. Then the lexicographic semigroupoids $R_\mu$ and
$R_\nu$ are isomorphic if and only if there is an order bijection $\pi :
\Omega \rightarrow \Omega^\prime$ such that $\nu (\pi (w)) = \mu (w)$
for all $w$ in $\Omega$.
\rm
\vspace{.3in}

\begin{mproof} Assume first that $\Omega$ and $\Omega^\prime$ have
no minimal elements.

Observe that a point $x = (x_w)$ in $X_\mu$ has a closed
orbit $\overline{\co (x)}$ which is a proper subset of $X_\mu$ if and
only if there exists an index of $t$ in $\Omega$ such that $x_w =1$ for
all $w \leq t$. Indeed, assume that  this does not hold and
consider an arbitrary point $y$ in $X_\mu$. For
an index $s \in \Omega$, with $x_s > 1$, let $z^s = (z_w)$ where $z_w =
y_w$ for $w < s, z_s = x_s -1$, and $z_w = x_w$ for $w > s$. If $A$ is
a basic open neighbourhood of $y$ then, from the assumption, it follows
that there exists an element $z^s$ in $A$. Since $(z^s, x) \in R_\mu$
the observation follows.

Divide the set of points with proper closed orbit into two types. A
point $X$ is of type 1 if there exists a first index $w_1$ such that
$x_{w_1} \neq 1$, and is of type 2 otherwise. We claim that $x$ is type
1 if and only if
$\overline{\co (x)}$ contains a closed orbit
$\overline{\co (y)}$ which is a relatively open-closed
proper subset. This property identifies the indexing in order topological
terms and is the basis of the proof.

If $x$ is a type 2 point then
there is an order interval decompositions $\Omega = \Omega_1 +
\Omega_2$, where $\Omega_2$ has no first element, such that $y$ belongs
to the closed orbit of $x$ if and only if $y = (y_w)$ with $y_w =1$ for $w \in
\Omega_1$.
To see this let $\Omega_1$ be the largest order ideal in $\Omega$ of
indices $w$ with $x_w = 1$. Thus $\Omega = \Omega_1 +
\Omega_2$, where $\Omega_2$ has no first element and for any index
$u$ in $\Omega_2$ there is a smaller index $v$ with $x_v > 1$.
In particular $\overline{\co (x)}$ contains the set

\[
(..1..\mid x_v-1\mid ..*..).
\]

The union of all such sets is dense in the set

\[
\overline{\co (x)} = (..1..\mid..*..)
\]

associated with the decomposition of $\Omega$. It follows that
$\overline{\co (x)} $ is precisely this latter set, as required.

A relatively open-closed subset of $\overline{\co (x)}$ contains
a basic relatively open-closed neighbourhood of the form

\[
(..1..|..*..|z_{w_1}|..*..|z_{w_2}|..*..,\dots ,..*..|z_{w_n}|..*..)
\]

\noindent with $w_1 < w_2 < \dots < w_n$. Such a set contains points $z =
(z_w)$
with $z_t = 2$ for all $t$ in $\Omega_2$ with $t < w_1$. Since
$\overline{\co (z}) = \overline{\co (x)}$ for such a point it follows that
$\overline{\co (x)}$ cannot contain, properly, a relatively open-closed
closed orbit.

We now identify the integers $\mu(w)$ in order-topological terms.

Let $E \subseteq X_\mu$ be an intersection of closed orbits of type 2
points which is not
itself a closed orbit of a type 2 point. Then, in view of the
description above of the closed orbits of type 2 points, E  has the form

\[
E^{(q)} = \bigcup^{\mu(q)}_{i=1} (..1..|i|..*..),
\]

\noindent for some $q$ in $\Omega$.
{}From this set and the relation $R_\mu$ we can discover $\mu(q)$
in order-topological terms as follows. There is a unique
$R_\mu$-invariant probability measure, $\lambda$ say,
on the set $E^{(q)}$. The sets

\[
E^{(q)}_k = \bigcup^k_{i=1} (..1..|i|..*..)
\]

\noindent are the only closed orbits contained in $E^{(q)}$ with
positive $\lambda $ measure, and the reciprocal of the measure of the smallest
such set is $\mu(q)$.

Suppose now that
$\alpha^{(2)} : R_\mu \rightarrow R_\nu$
is a semigroupoid isomorphism implemented by
$\alpha : X_\mu \rightarrow X_\nu$.
In particular $\ga$ is a homeomorphism.
If $p \in \Omega^\prime$ then write $E^{(p)}$ for the sets in $X_\nu$
that are analogous to the subsets
$E^{(q)}$.  Then
from the paragraph above it follows that
$\ga(E^{(q)}) = E^{\pi(q)}$ for some element
$\pi(q)$ in $\Omega^\pr$. The map  $\pi$ is an order isomorphism
and from the characterisation of the numbers $\mu(q)$ above
it follows that $\nu(\pi(q)) = \mu(q)$ for all $q$ in $\Omega$.

Suppose finally that $\Omega$ has a minimal element $w_0$.
Then all points of $X_\mu$ except the unique maximal
point $x_{max} = (\mu(w))$ have proper closed orbits.
Also, the converse is true. Note that there exist
open-closed orbits and a smallest open-closed orbit,
namely

\[
(1|..*..).
\]

Since this set is determined in order-topological terms we can restrict
considerations to this subset of $X_\mu$ and deduce the theorem
in this case from
the one already considered.

\end{mproof}
\newpage

\section{ The Case $\Omega = \omega^* + \omega$}

\noindent A classification is now given for the lexicographic semigroupoids
$R_\mu$ in the case where, in the standard notation,
$\Omega$ has order type $ \omega^* + \omega$.
This order type is quite a bit more subtle than that of the
dense orders in that semigroupoid isomorphisms may have to be effected
by homeomorphisms defined in terms of a recoding.

Throughout this section  assume that $\Omega = \ZZ \backslash \{0\}$
and that $\mu,\  \nu$ are maps from $\Omega$ to $\ZZ_+$ with

\[
\mu (k) = r_k,\  \mu (-k) = s_k,\  \nu (k) = t_k, \ \nu (-k) = u_k,
\]

\noindent for $k = 1, 2, \dots$. Associate with $\mu$ the pair $(\underline{r},
\underline{s})$ of generalised integers

\[\underline{r} = r_1 r_2 \dots,\ \
\underline{s} = s_1 s_2 \dots.
\]

\noindent  Define an equivalence relation $\sim$
on pairs of generalised (or finite) integers by $(\underline{r},
\underline{s}) \sim (\underline{t}, \underline{u})$ if and only if
$\underline{r}\underline{s} = \underline{t}\underline{u}$ and
there exist coprime natural numbers $a, b$ such that $b \underline{r} =
a \underline{t}$ and $a \underline{s} = b \underline{u}$.
\vspace{.3in}

\noindent {\bf Theorem 3} \it Let $\mu, \nu : \Omega \rightarrow \IN$, with
associated
pairs of generalised integers (possibly finite) $(\underline{r},
\underline{s})$ and $ (\underline{t}, \underline{u})$ respectively. Then
the lexicographic semigroupoids $R_\mu$ and $R_\nu$ are isomorphic if
and only if $(\underline{r}, \underline{s}) \sim (\underline{t},
\underline{u})$.
\vspace{.3in}
\rm

\begin{mproof}
We prove the necessity of the condition for isomorphism.
The sufficiency direction is
relatively straightforward
and is left to the reader. (See also \cite{hop-scp} and \cite{poo}.)

Let $\alpha : X_\mu \rightarrow X_\nu$ be a bijection such that
$\alpha^{(2)} : R_\mu \rightarrow R_\nu$ is a topological isomorphism. The
inversion map $\theta : X_\mu \times X_\mu \rightarrow X_\mu \times
X_\mu$ given by $\theta ((x, x^\prime)) = (x^\prime , x)$ is an
automorphism of the equivalence relation $\tilde{R}_\mu$, and
$\tilde{R}_\mu = R_\mu \cup \theta (R_\mu)$. It follows that
$\alpha^{(2)}$ maps $\tilde{R}_\mu$ homeomorphically onto
$\tilde{R}_\nu$. In particular, by the discussion in section 1, it
follows that the generalised integer for $\tilde{R}_\mu$ coincides with
that for $\tilde{R}_\nu$. That is, $\underline{r}\  \underline{s} =
\underline{t} \ \underline{u}$.

Let $X_{\mu, 0} \subseteq X_\mu$ be the set of points $x = (x_k)$
with $x_k =1$ for sufficiently small $k$.
If $X_{\mu ,0}$ is a proper
subset of $X_\mu$ then it is precisely the set of points in $X_\mu$
with proper closed orbits, and so $\alpha (X_{\mu , 0}) = X_{\nu ,0}$.
The set
$X_{\mu, 0} $ contains the special points for which, in addition,
 $x_k = r_k$ for all sufficiently large $k$. These special points
can be characterised
order topologically.
Indeed they are precisely the points $x$ for which there exists a
point $x^+$ whose closed orbit $ \overline{\co(x^+)}$ is the union of
$\{x+\}$ and $ \overline{\co(x)}$.
These so-called gap points are discussed in
\cite{hop-scp} and  \cite{scp-book}.
Thus, if $x_* = (\dots, 1, 1, \hat{r}_1, r_2, \dots)$, where
the symbol $\ \hat{}$\
\ \ indicates the coordinate position $k =1$, then $\alpha (x_*)$ may be
written as

\[
(\dots, 1, 1, w_{-j +1}, w_{-j}, \dots , w_{j-1}, u_j, u_{j+1}, \dots)
\]

\noindent for some positive integer $j$. We have

\[
\overline{\co (x_*)} = (\dots 1, 1, \hat{*}, *, \dots),
\]
\noindent and
\[
\overline{\co (\alpha (x_*))} = \{(\dots 1, w^\prime, y_j, y_{j+1},
\dots): w^\prime \in [1_j , w]\}
\]

\noindent where $y_k \leq u_k$ for $k \geq j$, and where $w^\prime$ is any word
of
length $2j -2$ which preceeds, or is equal to, the word $w = (w_{-j
+1}, \dots , w_{j-1})$. That is, $w^\prime$ belongs to the lexicographic
order interval $[1_j, w]$ where $1_j$ is the word with $2j-2$
coordinates all equal to 1. Let $n$ be the number of words in this order
interval and note that the restricted topological equivalence relation
$\tilde{R}_\nu \mid \overline{\co (\alpha (x_*))}$ has generalised
integer $nt_j t_{j+1} \dots$. Since $\alpha$ induces an isomorphism
between this relation and the topological equivalence relation
$\tilde{R}_\mu \mid \overline{\co (x_*)}$ then, once again, by the
discussion in section 1, we have $nt_j t_{j+1} \dots = r_1 r_2 \dots$,
and hence $m \underline{r} = n \underline{t}$ where $m = t_1 t_2 \dots
t_{j-1}$. Although $\underline{r} \underline{s} = \underline{t}
\underline{u}$ we cannot yet conclude that $(\underline{r},
\underline{s}) \sim (\underline{t}, \underline{u})$.

Define $d_\mu : X_{\mu, 0} \rightarrow \IR$ by

\[
d_\mu (x_k) = \sum^\infty_{k=1} \frac{x_k -1}{r_1 r_2 \dots r_k} +
\sum^\infty_{k=1} (x_{-k} -1) s_0 s_1 \dots s_{k-1},
\]

\noindent where $s_0 =1$, and define $d_\nu : X_{\nu, 0} \rightarrow \IR$
similarly. If $x \in \overline{\co (x_*)}$ then we can interpret
$d_\mu (x)$ as the
measure of $\overline{\co (x)}$ with respect to the unique normalised
$R_\mu-$invariant Borel measure on $\overline{\co (x_*)}$
Call this measure
$\lambda_\mu$, and note that there is a unique $\tilde{R}_\mu$-invariant
extension to $X_{\mu,0}$ which we also denote by $\lambda_\mu$.

Because of the uniqueness of normalised invariant measures it follows
that $\lambda_\nu \circ \alpha  = c \lambda_\mu$ and $d_\nu (\alpha (x))
= cd_\mu (x)$ for some positive constant $c$. Since $d_\mu (x_*) =1$ we
have $c = d_\nu (\alpha (x_*))$.
To see that $c = n/m$ note that there are precisely $m$ points in the
product $[t_1] \times \dots \times [t_{j-1}]$. Thus, from the definition
of $n$ and the product measure $\lambda_\nu$ we see that
 $\lambda_\nu(\overline{\co (\alpha (x_*))}) = n/m$.

We now use the connection $d_\nu (\alpha (x)) = cd_\nu (x)$  to show
that $n \underline{s} = m \underline{u}$. We can assume that
$ \underline{s}$ and $ \underline{u}$ are not finite.

Let $n/m = a/b$ where $(a, b) =1$ and suppose, by way of contradiction,
that $a \underline{s}$ does not divide $b \underline{u}$. Then, there
exists a prime number $p$ and a positive integer $g$ such that

\[
p^g \mid a\underline{s},\ \  p^{g-1} \mid b\underline{u}, \ \
\hbox{ and } p^g \not{\mid}
\ b\underline{u}.
\]

\noindent Choose $l$ large enough so that $as_1\dots s_l = k_1p^g$ for some
integer
$k_1$. Then

\[
k_1 \frac{b}{a} \frac{u_1 u_2 \dots u_l}{s_1 s_2 \dots s_l} = b
\frac{u_1^\prime u^\prime_2 \dots u^\prime_l}{p}
\]

\noindent and

\[
(p, bu^\prime_1 u^\prime_2 \dots u^\prime_l u_{l+1} \dots u_v) =1
\]

\noindent for all $v > l$, where $u^\prime_1, \dots u_l^\prime$ are factors of
$u_1, \dots u_l$ respectively. Note that by increasing $l$, if
necessary, and compensating with a multiple of $k_1$, we can arrange
that the product $s = s_1 s_2 \dots s_l$ satisfies $p^{-1} <1 -s^{-1}$.
These numerical relations lead to the contradiction that $\alpha (E) =
X_\nu$ where $E$ is the proper open-closed subset

\[
E = \bigcup_{w \in W} (\dots * \dots, w, \dots * \dots)
\]

\noindent where $W$ is the set of words
$(s_l^\prime , \dots , s^\prime_1)$ which
are strictly less than $(s_l, s_{l-1}, \dots, s_1)$ in the lexicographic
ordering.

To see this let $y$ be an arbitrary point of $X_\nu$ and let
$F_v (y)$ be the closed set

\[
F_v (y) = \{y^\prime \in X_\nu : y^\prime = (y^\prime_k), y^\prime_k = y_k
\hbox{ for } k \geq -v\}
\]

\noindent for $v = 1, 2, \dots$. The range of $d_\nu$ on $F_v (y) \cap X_{\nu,
0}$ is an arithmetic progression, namely,

\[
d_\nu (F_v (y) \cap X_{\nu, 0}) = \{ku_1 u_2 \dots u_v + \xi : k \in
\ZZ_+\}
\]

\noindent where $\xi = d_\nu (y^+)$, and where $y^+$ agress with $y$ in
coordinates
indexed by $\ZZ_+$ and is equal to 1 in the remaining coordinates.

On the other hand we have

\[
d_\mu (E \cap X_{\nu , 0}) = \bigcup^\infty_{k=1} [ks, ks + (s-1)]
\]

\noindent and so

\[
d_\nu (\alpha (E) \cap X_{\mu , 0}) = \bigcup^\infty_{k=1} [cks, cks +
c(s-1)].
\]

\noindent In particular, since $c = a/b$,
the set $d_\nu (F_v (y) \cap X_{\nu, 0})$
meets $d_\nu
(\alpha (E) \cap X_{\nu, 0})$ if and only if the sets $B_v$ and $A$
meet, where

\[
A = \bigcup^\infty_{k=1} [k, k+ (1-s^{-1})]
\]

\noindent and

\[
B_v = \{k \frac{b}{a} \frac{u_1 u_2 \dots u_v}{s} + \frac{\xi}{s} : k \in
\ZZ_+\}.
\]

By our earlier remarks, and our choice of $l$, if $v > l$ then $B_v$
contains the set

\[
B^\prime_v = \{k \frac{b u^\prime_1  u^\prime_2 \dots  u^\prime_l
u_{l+1} \dots u_v}{p} + \frac{\xi}{s} : k \in \ZZ_+\}.
\]

\noindent By the coprimality of $p$ and $b  u^\prime_1  u^\prime_2 \dots u_v$
this
set  contains certain  positive integral translates of the points

\[
\frac{u}{p} + \frac{\xi}{s} , \hbox{ for } i = 1, \dots, p.
\]

\noindent Since $p^{-1} < (1 -s^{-1})$ it follows that $B_v^\prime$ meets the
set
$A$. We have thus shown that

\[
d_\nu (F_v (y) \cap X_{\nu, 0}) \cap d_\nu (\alpha (E) \cap X_{\nu, 0})
\neq \emptyset
\]

\noindent and hence, by the openness of $\alpha (E)$,  that

\[
F_v (y) \cap \alpha (E) \neq \emptyset, \hbox{ for } v = 1, 2, \dots .
\]

\noindent Since $y$ is the unique point in the intersection of the sets $F_v
(y)$,
it follows that $y \in \alpha (E)$ and hence that $\alpha (E) = X_\nu$,
the desired contradiction.

We have shown that $a \underline{s}$ divides $b \underline{u}$. Since
$\alpha^{-1}$ is also an automorphism we conclude that $b \underline{u}$
divides  $a \underline{s}$ and hence that $b \underline{u} =a
\underline{s}$, and $(\underline{r}, \underline{s}) \sim (\underline{t},
\underline{u})$.
\end{mproof}
\vspace{.3in}

\noindent {\bf Remarks 1.}  Note that from the proof above it follows that the
restrictions $R_\mu \mid X_{\mu, 0}$ and $R_\nu \mid X_{\nu ,0}$ are
isomorphic topological binary relatations if and only if $R_\mu$ and
$R_\nu$ are isomorphic topological binary relations. We use this fact in
the next section.

\noindent {\bf 2.} As we mention later, there is a close
association between approximately finite
topological binary relations and  approximately finite
triangular operator algebras. The algebras associated with the
$\ZZ$-ordered relations above correspond to
the so-called
alternation algebras
considered by Hopenwasser and Power \cite{hop-scp} and
by Poon \cite{poo}.
Their classification is also given, as part of a more general
study, by Donsig and Hopenwasser \cite{don-hop}.
The argument we give above follows very closely
the method of \cite{scp-outer}, which was restricted to the case
$\mu = \nu$. (There is an inadequacy in the arithmetic progression
argument of \cite{scp-outer} which is corrected in the
somewhat more general argument above.)

\noindent {\bf 3.} In principle it should be possible
to reformulate the arguments of  Poon \cite{poo}
and Hopenwasser and Donsig \cite{don-hop} in terms of binary relations
to give alternative proofs of
Theorem 3.
The arguments we have given are also suited to other
situations and in particular to
the generalised alternation algebras
associated with Markov chains and subshifts. (See \cite{scp-book}.)
We intend to report more fully on this elsewhere.
However, let us note the following example from \cite{scp-book}.

Let $R_\mu$ be the lexicographic topological binary
relation for $\Omega = \ZZ$ with the function $\mu(w) = 2$ for all $w$.
Let $R_1$ (resp.  $R_2$) be the topological subrelation
defined on the symbol subspace $X_1$ (resp. $X_2$) of $X_\mu$ for which
the pair \ $00$ \ (resp. $11$) is forbidden. Then, in contrast to their
generated equivalence relations, $R_1$ and $R_2$
are not isomorphic. Similarly the one-sided subrelations are not isomorphic.

\vspace{.3in}

\noindent {\bf  Automorphisms}
\vspace{.3in}

Fix $\mu : \Omega \rightarrow \IN$ as above
with lexicographic semigroupoid $R_\mu$ and generalised integer pair
$(\underline{r},\underline{s})$.
Let $d$ be the number
of primes $p$ such that both $\underline{r}$ and $\underline{s}$
are divisible by $p^\infty$.
\vspace{.3in}

\noindent {\bf Theorem 4} \cite{scp-outer}
\it The semigroupoid automorphism group
Aut$(R_\mu)$ is isomorphic to the
restricted direct product $\ZZ^d$.
\rm
\vspace{.3in}

\begin{mproof}An element $\alpha^{(2)} \in Aut (R_\mu)$ is a
topological isomorphism
induced by a bijection $\alpha : X_\mu \rightarrow X_\mu$.
By the argument in
the proof of Theorem 3, specialized to the case $\mu = \nu$, we have

\[
d_\mu (\alpha (x)) = c d_\mu (x)
\]

\noindent for $x \in X_{\mu , 0}$, where $c = a/b$ with $(a, b) = 1$, $a
\underline{r} =  b \underline{r}$, and $b \underline{s} = a
\underline{s}$. Furthermore, if $c =1$ then $\alpha$ is trivial. These
conditions imply that $a^\infty$ and $b^\infty$ divide $\underline{r}$
and $\underline{s}$.
It follows that  $c$ has the form $p_1^{e_1} p_2^{e_2} \dots p_k^{e_k}$
where $e_i \in \ZZ$ and each $p_i$ divides $\underline{r}$ and
$\underline{s}$ with infinite multiplicity. The mapping $\alpha
\rightarrow c$ gives the desired isomorphism.
\end{mproof}

Suppose that $p$ is a prime such that
$p^\infty$ divides $\underline{r}$ and $\underline{s}$. We note one
way in which the order-preserving homeomorphism $\alpha$ corresponding
to $p^{-1}$ may be identified.

Write $\underline{r} = pr_1 pr_2 pr_3 \dots, $ and
$\underline{s} = s_1 p s_2 p
s_3  \dots$ and obtain the identification

\[ \begin{array}{llc}
X_\mu & = & ( \dots [s_2] \times [p] \times [s_1]) \times ([p] \times [r_1]
\times [p] \times \dots)\\
& = & X_\mu^- \times X_\mu^+.
\end{array} \]

\noindent This is naturally isomorphic to

\[ \begin{array}{llc}
X_\lambda & = & ( \dots [p] \times [s_1] \times [p]) \times ([r_1]
\times [p] \times [r_2] \times \dots )\\
& = & X_\lambda^1 \times X_\lambda^+
\end{array} \]

\noindent by a map $\beta = \beta^- \times \beta^+$ which respects the factors
and
which induces the natural semigroupoid isomorphism $R_\mu \rightarrow
R_\mu$.  Let $\gamma : X_\lambda \rightarrow X_\mu$ be the right shift
homeomorphism. Then $\gamma \circ \beta$ is an automorphism of $X_\mu$
and its associated constant is $p^{-1}$.

\section{\bf  Countable Linear Orderings}

\noindent Let $\Omega$ be a countable linear ordering and define an equivalence
relation $\approx$ on $\Omega$ such that $w \approx v$ if the order
intervals $[w, v]$ and $[v, w]$ are finite. Then the set $\Omega /
\approx$ of equivalence classes is linearly ordered and each
equivalence class $< x>$
is itself a linearly ordered set which is isomorphic to a finite set or to one
of
$\ZZ_+, \ZZ_{-}$ and $\ZZ$.
Let $\mu : \Omega \to \{2,3,\dots\}$. Then to each class $<x>$ we can associate
a pair $p_\mu(<x>) = (\underline{r}, \underline{s})$, as in
section 3, consisting of finite or generalized integers.
\vspace{.3in}

\noindent {\bf Theorem 5}. \  \it Let $\Omega,\  \Lambda$ be
countable linear orderings with maps
$\mu : \Omega \rightarrow \{2, 3, \dots \},  \nu : \Lambda \rightarrow
\{2, 3, \dots \}$. Then the lexicographic semigroupoids $R_\mu$ and
$R_\nu$ are isomorphic if and only if there is an order preserving
bijection $\pi : \Omega / \approx \ \rightarrow \ \Lambda / \approx$ \ such
that $p_\mu (\pi (<w>)) \sim p_\mu (<w>)$ for all classes $<w>$ in
$\Omega / \approx$.
\rm
\vspace{.3in}

\begin{mproof}  As in the proof of Theorem 2  declare a point $x = (x_w)$
in $X_\mu$ to be a type 1 point if there exists a first coordinate $w$
for which $x_w \neq 1$. Once again, as in the proof of Theorem 2, these
points are identifiable in order-topological terms as those for which
the closed (half) orbit $\overline{\co (x)}$ properly contains relatively
open-closed closed orbits. If $w_2 > w_1$ then say that $w_1$ and
$w_2$ are {\it finitely equivalent} if there exist type 1 points $x, y$
with first non-unit coordinates $x_{w_1}$
and $y_{w_2}$ respectively such that $\overline{\co (y)}$ has positive measure
with respect to the unique invariant probability Borel measure
on $\overline{\co (x)}$.

Note that $w_1$ and $w_2$ are finitely equivalent if and only if
$w_1 \approx w_2$.
In view of this it follows that if $\alpha^{(2)} :
R_\mu \rightarrow R_\nu$ is a topological isomorphism induced by the
homeomorphism $\alpha$ then $\alpha$ induces a map $\pi : \Omega /
\approx \ \rightarrow \ \Lambda / \approx$.
Indeed, if $x$ is a type 1 point associated with $w$ in $\Omega$
then $\alpha(x)$ is a type one point associated with $u$ in $\Lambda$ and
we may define $\pi(<w>) = <u>$. Since the equivalence relation
$\ \approx$\  coincides with finite equivalence this is a well-defined
bijection.

Fix a class $<w>, $ for some $w \in \Omega$, and define the set

\[
X_\mu (<w>) = \bigcup_{x \in T(w)} \overline{\co (x)}
\]

\noindent where the union is taken over the set $T(w)$ of type 1 points $x$
associated with the class $<w>$. Then $\alpha$ restricts to a
homeomorphism

\[
\beta : X_\mu (<w>) \rightarrow X_\nu (\pi (<w>)).
\]

\noindent By considering the restriction of $\alpha^{(2)}$ to the set

\[
R_\mu \cap
(X_\mu (<w>) \times X_\mu  (<w>))
\]

\noindent we shall show that the lexicographic
semigroupoid $R(<w>, \mu)$ associated with the linear order $<w>$ and
the function $\mu$ is isomorphic to $R (\pi (<w>),\nu)$.
In view of the result in
section 3 this will complete the proof.

To this end define the
equivalence relation $E_\mu$ on $X_\mu (<w>)$ as the set of pairs
$(x,y)$ for which $x_u = y_u$ for all $u$ in $<w>$, and similarly define
$E_\nu$ on $X_\nu (\pi  (<w>))$. Thus, the set of equivalence classes
$X_\mu  (<w>) /E_\mu$ is isomorphic to the set

\[
X_{<w>} = \prod_{d \in <w>} [\mu (d)].
\]

\noindent Furthermore $R_\mu$ and $E_\mu$ induce a binary relation, $R_\mu /
E_\mu$ say, on the set $X_\mu (<w>)/E_\mu. $ That is
$ (x,y) \in R_\mu / E_\mu$
if and only if there exist $x^\prime, y^\prime$ in $X_\mu (<w>)$ with
$(x^\prime , y^\prime) \in R_\mu, x^\prime \in x$ and $y^\prime \in y$.
Also, under the natural identification above $R_\mu / E_\mu$ is the
lexicographic semigroupoid $R (<w>, \mu )$ on
$X_{<w>}$.

In view of these identifications it will be enough to show that the
equivalence relation $E_\mu$ can be defined in an order-topological
fashion and that the restriction of $\alpha$ to $X_\mu  (<w>) $ induces
a semigroupoid isomorphism from $R (<w>, \mu)$ to $R ( \pi  (<w>) ,
\mu)$.

In section 2 we saw that, up to a constant multiplier, the set $X_\mu
(<w>) $ carries a unique $R_\mu$-invariant measure, $\lambda_\mu$ say. If
$x, y \in X_\mu  (<w>) $ and $\lambda_\mu (\overline{\co (x)}) =
\lambda_\mu (\overline{\co (y)})$ then there are two possibilities;
either $x_u = y_u$ for all $u$ in $<w>$, or $x_u$ and $y_u$ correspond
to ``rational points'' in the sense that the symmetric
difference \ \ $\overline{\co
(x)} \Delta \overline{\co (y)}$ \ \ is a singleton, namely $\{x\}$ or
$\{y\}$. It follows that $E_\mu$ can be defined purely in
order-topological terms and hence that $\alpha$ maps the $E_\mu$-equivalence
classes to $E_\nu$-equivalence classes.

Since the given map $\alpha : X_\mu \rightarrow X_\nu$ is continuous and
since the basic open-closed sets of $X_\mu$ and $X_\mu$ generate the topology
of $X_\mu$ and $X_\mu$, respectively, it follows that the induced map

\[
\alpha_i : X_\mu (<w>) / E_\mu \rightarrow X_\nu (\pi (<w>)) / E_\nu
\]

\noindent is
bicontinuous. Indeed $\alpha$ maps an $E_\mu$-saturated basic open-closed set
to an $E_\nu$-saturated open-closed set and this is necessarily a finite union
of $E_\nu$-saturated basic open-closed sets. Similarly, it follows that
$\alpha_i^{(2)}$ induces a semigroupoid isomorphism from $R(<w>, \mu)$ to
$R(\pi (<w>), \mu)$.
\end{mproof}

The arguments given above and in the previous section are also effective in the
setting of
infinite lexicographic products of partially ordered
sets. This is illustrated in the following theorem in the case $\Omega = \IQ$
for which the proof in section 2 is applicable with little change.

For each rational $q \in \IQ$ let $\le_q$ be a connected partial
ordering on the finite set
$\{1,\dots,\mu(q)\}$. Then the product space $X_\mu$ carries
a natural semigroupoid

\[
R =
R(\{ \le_q : q \in \IQ\})
\]

\noindent which is the subset
of $R_\mu$ associated with the given
partial orderings. That is, $(x,y) \in R$
if and only if $(x,y) \in \tilde{R}_\mu$ and
$x_{w_0} \le_q y_{w_0}$\ for the smallest index $w_0$
such that $x_w \ne y_w.$

 \vspace{.3in}

\noindent {\bf Theorem 6} \ \it The semigroupoids
$R(\{ \le_q : q \in \IQ\})$ and
$R(\{ \prec_q : q \in \IQ\})$ are isomorphic if and only if there is an order
bijection $\pi : \IQ \to \IQ$ such that the partial orderings
$\le_q$ and $\prec_{\pi(q)}$ are isomorphic for all $q \in \IQ$.
\rm

\vspace{.3in}

\section{\bf Applications to Operator Algebras}

The operator algebra $T_n$ is the subalgebra of the complex matrix algebra
$ M_n$ consisting of upper triangular complex matrices
and endowed with  the usual operator norm. These algebras $A$
are {\it triangular} in the sense that $A \cap A^*$ is a maximal abelian
self-adjoint subalgebra. Recently there has been considerable interest
in classifying the many diverse families of triangular operator
algebras arising as direct limits of these finite-dimensional algebras
and their direct sums. (See \cite{scp-book}.)

The following construction is given in \cite{scp-lex1}.

Let $\Omega, \mu$ be as above, with $n_w = \mu (w)$. Let $F \subseteq \Omega$
be a finite
subset, say $w_1 < w_2 < \dots < w_k$, and let $w_t < w < w_{t+1}$, for some
$t$. Set $G = F \cup \{w\}$, $n_F = n_{w_1} n_{w_2}\dots n_{w_k}$, and $n_G =
n_\omega n_F$.   Define a unital algebra injection $\phi_{F,G} : T_{n_F} \to
T_{n_G}$ as follows. View $T_{n_F}$ as the (maximal triangular) subalgebra
of $M_{n_{w_1}} \otimes \dots \otimes M_{n_{w_k}}$ which is spanned by the
matrix units

\[
e_{{{\bf i}},{{\bf j}}} = e_{i_1, j_1} \otimes \dots \otimes e_{i_k,j_k}
\]

where the multi-index ${{\bf i} } = (i_1, \dots, i_k)$ precedes ${{\bf j}} =
(j_1, \dots,
j_k)$ in the lexicographic ordering. Thus either ${{\bf i}} =  {{\bf j}}$ or
the first
$i_p$ differing from $j_p$ is strictly less than $j_p$. Similarly  identify
$T_{n_G}$ for the ordered subset $G$ and  set $\phi_{F,G}$ to be the linear
extension of the correspondence

\[
e_{{{\bf i}},{{\bf j}}} \to \sum^{n_\omega}_{s=1} e_{i_1, j_1} \otimes \dots
\otimes
e_{i_t, j_t} \otimes e_{s,s} \otimes e_{i_{t+1},j_{t+1}} \otimes \dots \otimes
e_{i_{k},j_{k}}
\]

In a similar way, or by composing maps of the above type, define
$\phi_{F,G}$ for $F \subseteq G$, general finite subsets. These maps are
isometric and so determine the Banach  algebra

\[
A(\Omega, \nu) = {\lim_{\to}}_{F \in {\cal F}} T_{n_F}
\]

\noindent where the direct limit is taken over the directed set ${\cal F}$ of
finite
subsets of $\Omega$. Each $\phi_{F,G}$ has an extension to a
C*-algebra injection from $M_{n_F}$ to $M_{n_G}$ and so
it follows that we may
view $A(\Omega, \nu)$ as a closed unital subalgebra of the UHF
C*-algebra $B(\Omega, \nu) = {\displaystyle \lim_\to M_{n_F}}$.
\vspace{.3in}

\noindent {\bf Theorem 7} \  \it The following statements are equivalent.

(i) $A(\Omega,\mu)$ and $A(\Lambda,\nu)$ are isometrically isomorphic Banach
algebras.

(ii) $R_\mu$ and $R_\nu$ are isomorphic lexicographic semigroupoids.

(iii) There is an order preserving bijection
 $\pi : \Omega / \approx \ \rightarrow \ \Lambda / \approx$ \ such
that
\[
p_\mu (\pi (<w>)) \sim p_\mu (<w>)
\]
for all classes $<w>$ in
$\Omega / \approx$.

\rm
\vspace{.3in}

\begin{mproof} The semigroupoids are readily identifiable with the
topological fundamental relations of the operator algebras. (See, for example,
the discussions of \cite{hop-scp} and \cite{scp-book}.) The equivalence of
(i) and (ii) is now immediate from Theorem 7.5 of \cite{scp-book} and so
Theorem 5 completes the equivalences.
\end{mproof}

An immediate corollary of the last theorem is
that there are uncountably many triangular algebras $A$ of the form
$A(\Omega,\mu)$ with
$C^*(A)$ equal to the $2^\infty$ UHF C*-algebra.
On the other hand there are only three such (infinite-dimensional)
algebras of the form $A(\ZZ,\mu)$,
namely the pure refinement algebra, the standard limit algebra
and the alternation algebra with invariant $(2^\infty,2^\infty)$.

We now explain how the  algebras above can also be interpreted
in term of a lexicographic product operation at the algebraic level,
as described in \cite{scp-lex1}.

Let $A$ be an operator algebra admitting a subdiagonal
decomposition in the sense that

\[
A \ \ = \ \ A \cap A^* \ +\ A^0
\]

\noindent where $A \cap A^*$ is a maximal abelian subalgebra of $A$ ,
and $A^0$ is the kernel of a contractive
homomorphism $A \to A \cap A^*$. In particular,
this holds if $A$ is a regular triangular subalgebra
of an AF C*-algebra (\cite{scp-book} or, more generally, if
$A$ is a subdiagonal algebra in the sense of Arveson \cite{arv}.
If $A$ and $B$ are triangular operator algebras
admitting such decompositions then define their
{\it lexicographic product} $A \star B$ to be the closed subalgebra
of the injective tensor product
$C^*(A) \ot C^*(B)$ given by

\[
A \star B =  (A \cap A^*) \ot B \ +\  A^0 \ot C^*(B).
\]

One can verify that the inclusions $\phi_{F,G}$ defined above
coincide with the natural inclusions

\[
T_{n_{w_1}}\star \dots \star T_{n_{w_k}}
\ \ \to \ \ T_{n_{w_1}} \star \dots
T_{n_{w_t}}\star T_{n_{w}}\star T_{n_{w_{t+1}}}\star \dots \star T_{n_{w_k}}.
\]

\noindent In fact the lexicographic product is an associative operation
and the algebras $A(\Omega,\mu)$ can be viewed (unambiguously) as
infinite lexicographic products of upper triangular
matrix algebras over the ordering $\Omega$.
The following theorem is a corollary of Theorem 6.
\vspace{.3in}
\newpage

\noindent {\bf Theorem 8}\ \it Let $G_q$ and $H_q$ for $q$
in $\IQ$ be connected transitive antisymmetric
digraphs with triangular digraph algebras $A(G_q)$ and $A(H_q)$. Then the
lexicographic products

\[
\prod_{q \in \  \IQ} \star A(G_q)  \ \ \ \ \mbox{ and }\ \ \ \
  \prod_{q \in \ \IQ} \star A(H_q)
\]

\noindent are isometrically isomorphic triangular operator algebras if and
only if there is an
order bijection $\pi$ such that the digraphs
$G_q$ and $H_{\pi(q)}$ are isomorphic
for all rationals $q$.
\rm

\end{document}